\documentclass[a4paper,12pt]{article}
\pdfoutput=1
   \usepackage[top=2.5cm,bottom=2.5cm,left=2.5cm,right=2.5cm]{geometry}
\usepackage{cite,amsmath,amssymb}
\usepackage[margin=1cm,%
            font=small,%
            format=hang,%
            labelsep=period,%
             labelfont=bf]{caption}
\usepackage{graphicx}
\usepackage{caption}
\usepackage{float}
\usepackage{amsmath}
\usepackage{cases}
\usepackage{mathrsfs}
\usepackage{amsfonts}
\usepackage{color}
\usepackage{setspace}
\usepackage{cite}
\usepackage{amsmath}
\usepackage{array}
\usepackage{booktabs}

\usepackage{graphicx}
\usepackage{subfigure}
\usepackage{CJK}
\usepackage{caption}
\usepackage{float}
\usepackage{amsmath}
\usepackage{cases}
\usepackage{mathrsfs}
\usepackage{amsfonts}
\usepackage{times}
\usepackage{mathptmx}
\usepackage{cite}
\usepackage{setspace}
\usepackage{times}
\usepackage{mathptmx}
\usepackage{geometry}

\renewcommand{\figurename}

\renewcommand{\figurename}

\makeatletter
\def\thanks#1{\protected@xdef\@thanks{\@thanks
        \protect\footnotetext{#1}}}
\makeatother

\begin{document}
\title{\textbf{Arithmetic-Geometric Spectral Radius of Trees and Unicyclic Graphs}
\vspace*{0.3cm}}
\author{ \bf  Ruiling Zheng, \, Xian$^{^{\textbf{,}}}$an Jin\\
\small\em School of Mathematical Sciences, Xiamen University, Xiamen  361005, P. R. China\\ \small\em E-mail:
rlzheng2017@163.com, xajin@xmu.edu.cn.\vspace*{0.3cm}}
 \date{}
\maketitle\thispagestyle{empty}
\vspace*{-1cm}

\begin{abstract}
The arithmetic-geometric matrix $A_{ag}(G)$ of a graph $G$ is a square matrix, where the $(i,j)$-entry is equal to $
\displaystyle \frac{d_{i}+d_{j}}{2\sqrt{d_{i}d_{j}}}$ if the vertices $v_{i}$ and $v_{j}$ are adjacent, and 0 otherwise. The
arithmetic-geometric spectral radius of $G$, denoted by $\rho_{ag}(G)$, is the largest eigenvalue of the
arithmetic-geometric matrix $A_{ag}(G)$. Let $S_{n}$ be the star of order $n\geq3$ and $S_{n}+e$ be the unicyclic graph obtained from $S_{n}$ by adding an edge. In this paper, we prove that for any tree $T$ of order $n\geq2$,

\begin{center}
$\displaystyle 2\cos\frac{\pi}{n+1}\leq\rho_{ag}(P_{n})\leq\rho_{ag}(T)\leq\rho_{ag}(S_{n})=\frac{n}{2},$
\end{center}

\noindent with equality if and only if $T\cong P_{n}$ for the lower bound, and if and only if $T\cong S_{n}$ for the upper
bound. We also prove that for any unicyclic graph $G$ of order $n\geq3$,

\begin{center}
$\displaystyle 2=\rho_{ag}(C_{n})\leq\rho_{ag}(G)\leq\rho_{ag}(S_{n}+e),$
\end{center}

\noindent the lower (upper, respectively) bound is attained if and only if $T\cong C_{n}$ ($T\cong S_{n}+e$, respectively) and
$\displaystyle\rho_{ag}(S_{n}+e)<\frac{n}{2}$ for $n\geq7$.

\medskip

\end{abstract}

\vspace*{0.35cm}
\baselineskip=0.30in
\begin{spacing}{1.5}
\section{Introduction}

\ \ \ \ \ Throughout this paper we consider finite, undirected and simple graphs. Let $G$ be a graph with vertex set
$V(G)=\{v_{1},v_{2},\ldots,v_{n}\}$ and edge set $E(G)$. An edge $e\in E(G)$ with end vertices $v_{i}$ and $v_{j}$ is denoted by
$v_{i}v_{j}$. For $i=1,2,\ldots,n$, we denote by $d_{i}$ the degree of the vertex $v_{i}$ in $G$. As usual, let $P_{n}$, $S_{n}$ and
$C_{n}$ be the path, the star and the cycle of order $n\geq3$. Let $S_{n}+e$ with $n\geq3$ denote the unicyclic graph obtained
from $S_{n}$ by adding an edge. The spectral radius of a complex matrix $M$ is the largest value among the modules of all eigenvalues of $M$. In the case that $M$ is real symmetric, its spectral radius is exactly the largest eigenvalue and we denote it by $\displaystyle \rho(M)$. The spectral radius of the adjacency matrix
$A(G)=(a_{i,j})$ of $G$ is referred as the spectral radius of $G$.

In 1975, a so-called Randi\'{c} index was proposed\cite{1}, it is
defined as

\begin{center}
$ \displaystyle R_{-1}(G)=\sum\limits_{v_{i}v_{j}\in E(G)}\frac{1}{d_{i}d_{j}}.$
\end{center}

\noindent The Randi\'{c} matrix $R(G)=(r_{i,j})$ of $G$, where $\displaystyle
r_{i,j}=\frac{1}{d_{i}d_{j}}$ if $v_{i}v_{j}\in E(G)$ and 0 otherwise, seems to be
first time used in 2005 by Rodr\'{\i}guez, who referred to it as the ``weighted adjacency matrix" \cite{555} and
the ``degree adjacency matrix" \cite{666}.  Moreover, the role of Randi\'{c} matrix in the Laplacian theory was
clarified in \cite{77,999}. The Randi\'{c} spectral radius was studied in \cite{2,888,77,999} and the references cited therein.

In 1994, Yang et al.\cite{3} proposed the extended matrix of $G$, denoted by $A_{ex}(G)=(c_{i,j})$,  where
$\displaystyle c_{i,j}=\frac{1}{2}\Big(\frac{d_{i}}{d_{j}}+\frac{d_{j}}{d_{i}}\Big)$ if $v_{i}v_{j}\in E(G)$ and $0$ otherwise.
The corresponding topological index is the symmetric division deg index, which is proposed by Vuki\v{c}evi\'{c}\cite{4}
and
written as

\begin{center}
$ \displaystyle SDD(G)=\sum\limits_{v_{i}v_{j}\in E(G)}\Big(\frac{d_{i}}{d_{j}}+\frac{d_{j}}{d_{i}}\Big).$
\end{center}

\noindent  The extended spectral radius were studied in \cite{5,6,7} and the references cited therein.

In 1998, Estrada, Torres, Rodr\'{\i}guez and Gutman\cite{8} introduced the atom-bond connectivity $(ABC)$ index, defined as

\begin{equation*}
\displaystyle ABC(G)=\sum\limits_{v_{i}v_{j}\in E(G)}\sqrt{\frac{d_{i}+d_{j}-2}{d_{i}d_{j}}}.
\end{equation*}

\noindent The ABC matrix of G, put forward by Estrada\cite{10}, is defined to be $\Omega(G)=(\omega_{i,j})$,  where
$\displaystyle \omega_{i,j}=\sqrt{\frac{d_{i}+d_{j}-2}{d_{i}d_{j}}}$ if $v_{i}v_{j}\in E(G)$ and $0$ otherwise. Chen \cite{11} showed
that for
any tree $T$ of order $n\geq3$,

\begin{center}
$\displaystyle\sqrt{2}\cos\frac{\pi}{n+1}\leq\rho(\Omega(T))\leq\sqrt{n-2},$
\end{center}

\noindent with left (right, respectively) equality if and only if $T\cong P_{n}$( $T\cong S_{n}$, respectively), which was conjectured
in\cite{12}. Recently, Li and Wang\cite{13} proved that for an unicyclic graph $G$ of order $n\geq4$,

\begin{center}
$\displaystyle \sqrt{2}=\rho(\Omega(C_{n}))\leq\rho(\Omega(G))\leq\rho(\Omega(S_{n}+e)),$
\end{center}

\noindent with equality if and only if $T\cong C_{n}$ for the lower bound, and if and only if $T\cong S_{n}+e$ for the upper bound,
which was conjectured in\cite{14}. Yuan and Du characterized the bicyclic graphs with the largest and second largest $ABC$ spectral radius, which are just the unique two bicyclic graphs of order $n$ with maximum degree $n-1$ for $n\geq7$\cite{30}.

In 2015, Shegehall and Kanabur\cite{15} proposed the arithmetic-geometric index of $G$, which is defined as

\begin{center}
$ \displaystyle AG(G)=\sum\limits_{v_{i}v_{j}\in E(G)}\frac{d_{i}+d_{j}}{2\sqrt{d_{i}d_{j}}}.  $
\end{center}

\noindent   Recently, Zheng et al.\cite{16} considered the arithmetic-geometric matrix (AG matrix) of $G$, denoted by
$A_{ag}(G)=(h_{i,j})$,
where $\displaystyle h_{i,j}=\frac{d_{i}+d_{j}}{2\sqrt{d_{i}d_{j}}}$ if $v_{i}v_{j}\in E(G)$ and 0 otherwise. We denote $\rho(A_{ag}(G))$  by $\rho_{ag}(G)$, and call it
the arithmetic-geometric spectral radius ($AG$ spectral radius) of $G$. Moreover, Zheng et al.\cite{16} gave several bounds for
the AG spectral radius in terms of the maximum degree, the minimum degree of $G$, the Randi\'{c} index $ \displaystyle
R_{-1}(G)$ and the first Zagreb index $M_{1}(G)=\sum\limits_{v_{i}v_{j}\in E}(d_{i}+d_{j})$. Guo and Gao \cite{17} offered several other
bounds for the $AG$ spectral radius as well.

In this paper, we consider the $AG$ spectral radius of trees and unicyclic graphs.  We prove that:
(1) for any tree $T$ of order $n\geq2$,

\begin{center}
$\displaystyle 2\cos\frac{\pi}{n+1}\leq\rho_{ag}(P_{n})\leq\rho_{ag}(T)\leq\rho_{ag}(S_{n})=\frac{n}{2},$
\end{center}

\noindent with equality if and only if $T\cong P_{n}$ for the lower bound, and if and only if $T\cong S_{n}$ for the upper
bound.

\noindent (2) for any unicyclic graph $G$ of order $n\geq3$,

\begin{center}
$\displaystyle 2=\rho_{ag}(C_{n})\leq\rho_{ag}(G)\leq\rho_{ag}(S_{n}+e),$
\end{center}

\noindent the lower (upper, respectively) bound is attained if and only if $T\cong C_{n}$ ($T\cong S_{n}+e$, respectively) and
$\displaystyle\rho_{ag}(S_{n}+e)<\frac{n}{2}$ for $n\geq7$.

\section{Lemmas}
\noindent {\bf{Lemma 2.1}}\cite{19}
Let $M$ ba a nonnegative matrix of order $n$. Let $\textbf{x}$ be a positive column vector of dimension
$n$, i.e.
every entry is positive. If $k>0$ such that $M\textbf{x}\leq k\textbf{x}$, then $\rho(M)\leq k$.

\noindent {\bf{Lemma 2.2}}\cite{19}
Let $B=(b_{ij})$ and $D=(d_{ij})$ be two nonnegative matrices of order $n$. If $B\geq D$, i.e., $b_{ij} \geq d_{ij}\geq0$
for all $i,j$, then $\rho(B) \geq\rho(D)$.

\noindent {\bf{Lemma 2.3}}\cite{19}
For any principle submatrix $M^{\prime}$ of a nonnegative
matrix $M$ of order $n$, $\rho(M^{\prime})\leq\rho(M)$.

Let $DT(p,q)$ denote the double star obtained by joining the center of $K_{1,p}$ to that of $K_{1,q}$, which is depicted in Fig. 1. The
following Lemma was due to Chang and Huang\cite{21}.

\begin{figure}[H]
\begin{center}
    \includegraphics[width=8cm]{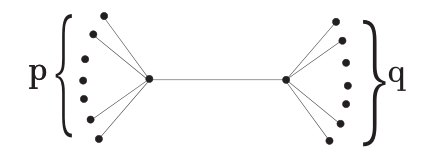} 
  \end{center}
  \caption*{\small \textbf{Fig. 1.} The double star $DT(p,q).$}
\end{figure}

\noindent {\bf{Lemma 2.4}}
For two integers $p,q$ with $p\geq q\geq 1$,

\begin{center}
$ \displaystyle\rho(A(DT(p,q)))=\sqrt{\frac{p+q+1+\sqrt{(p+q+1)^{2}-4pq}}{2}}.$
\end{center}

\noindent {\bf{Lemma 2.5}}\cite{22} If $G$ is an unicyclic graph, then
\begin{center}
$\displaystyle 2=\rho(A(C_{n}))\leq\rho(A(G))\leq\rho(A(S_{n}+e)),$
\end{center}

\noindent the lower bound is attained if and only if $G\cong C_{n}$, the upper bound is obtained if and only if $G\cong S_{n}+e$.

\noindent {\bf{Lemma 2.6}}\cite{20}
Let $\displaystyle \Phi(P_{n};\lambda)$ be the characteristic polynomial of $P_{n}$. Then

\begin{center}
$\displaystyle \Phi(P_{n};\lambda)=\frac{1}{\sqrt{\lambda^{2}-4}}\Big(x_{1}^{n+1}-x_{2}^{n+1}\Big)$,
\end{center}

\noindent where
$\displaystyle x_{1}=\frac{1}{2}\Big(\lambda+\sqrt{\lambda^{2}-4}\Big)$ and $\displaystyle
x_{2}=\frac{1}{2}\Big(\lambda-\sqrt{\lambda^{2}-4}\Big)$ are the roots of the equation $x^{2}-\lambda x+1=0.$

Zheng et al.\cite{16} gave the following upper bound for the AG spectral radius of graphs with order $n\geq2$ and size $m$.

\noindent {\bf{Lemma 2.7}}\cite{16}
Let $G$ be a graph of order $n\geq2$ and size $m$. Then

\begin{center}
$\displaystyle\rho_{ag}(G)\leq\frac{1}{2}\Big(\sqrt{n-1}+\frac{1}{\sqrt{n-1}}\Big)\sqrt{2m-n+1},$
\end{center}
 where the equality holds if and only if $G\cong S_{n}$.

 \section{The $AG$ spectral radius of trees}

 \noindent {\textbf{Theorem 3.1.}}
If $T$ is a tree of order $n\geq2$, then

\begin{center}
$\displaystyle 2\cos\frac{\pi}{n+1}<\rho_{ag}(P_{n})\leq\rho_{ag}(T)\leq\rho_{ag}(S_{n})=\frac{n}{2},$
\end{center}

\noindent with equality if and only if $T\cong P_{n}$ for the lower bound, and if and only if $T\cong S_{n}$ for the upper bound.

\noindent {\textbf{Proof.}} If $n=2,3$, then it is obvious. Now we suppose that $n\geq4$. According to Lemma 2.7, it follows that

\begin{center}
$\displaystyle\rho_{ag}(T)\leq\rho_{ag}(S_{n})=\frac{n}{2},$
\end{center}

\noindent we get the upper bound.

Then it suffices to show that for any tree $T\ncong P_{n}$,
\begin{center}
$\displaystyle2\cos\frac{\pi}{n+1}<\rho_{ag}(P_{n})< 2 \leq \rho_{ag}(T).$
\end{center}

We first prove that $\displaystyle2 \leq \rho_{ag}(T)$. If $T\ncong P_{n}$, then it has at least one vertex of degree greater than 2. We consider the following
two cases.

\noindent {\textbf{Case 1.}} $T$ has two adjacent vertices of degree greater than 2.

 Without loss of generality, we assume that $v_{p}v_{q}$ is an edge of $G$ such that $d_{p}\geq d_{q}\geq3$. Clearly,
 $DT:=DT(d_{p}-1,d_{q}-1)$ is an induced subgraph of $T$. Let $\displaystyle
 h(x,y)=\frac{x+y}{2\sqrt{xy}}$, then $\displaystyle h(x,y)\geq1$ for $x,y>0$. By Lemma 2.2, We get
 $\displaystyle\rho_{ag}(T)\geq\rho(A(T))$. And by Lemmas 2.3, 2.4, we have

 \begin{center}
$\displaystyle\rho_{ag}(T)\geq\rho(A(T))\geq\rho(A(DT))=\sqrt{\frac{d_{p}+d_{q}-1+\sqrt{(d_{p}+d_{q}-1)^{2}-4(d_{p}-1)(d_{q}-1)}}{2}}.$
\end{center}

\noindent Then we just need to prove that

\begin{center}
$\displaystyle \sqrt{\frac{d_{p}+d_{q}-1+\sqrt{(d_{p}+d_{q}-1)^{2}-4(d_{p}-1)(d_{q}-1)}}{2}}\geq2$ for $d_{p}\geq d_{q}\geq3$.
\end{center}

\noindent Since $d_{p}\geq d_{q}\geq3$, we obtain

\begin{center}
$\displaystyle d_{p}+d_{q}-1+\sqrt{(d_{p}+d_{q}-1)^{2}-4(d_{p}-1)(d_{q}-1)}$
\end{center}

$\ \ \ \ \ \ \ \  \ \ \ \ \ \ \ \ \ \ \ \ \ \ \ \ \ \ \geq$5+$\sqrt{(d_{p}+d_{q}-1)^{2}-4(d_{p}-1)(d_{q}-1)}$

$\ \ \ \ \ \ \ \  \ \ \ \ \ \ \ \ \ \ \ \ \ \ \ \ \ \ =5+\sqrt{(d_{p}-d_{q})^{2}+2d_{p}+2d_{q}-3}$

$\ \ \ \ \ \ \ \  \ \ \ \ \ \ \ \ \ \ \ \ \ \ \ \ \ \ \geq5+\sqrt{2d_{p}+2d_{q}-3}$

$\ \ \ \ \ \ \ \  \ \ \ \ \ \ \ \ \ \ \ \ \ \ \ \ \ \ \geq5+3=8$.

\noindent Hence

\begin{center}
$\displaystyle \sqrt{\frac{d_{p}+d_{q}-1+\sqrt{(d_{p}+d_{q}-1)^{2}-4(d_{p}-1)(d_{q}-1)}}{2}}\geq2$.
\end{center}

\noindent {\textbf{Case 2.}} $T$ has no two adjacent vertices of degree greater than 2.

 In this case, because $T\ncong P_{n}$, then $T$ has at least one vertex of degree greater than 2 and its adjacent vertices must be of
 degree 2 or 1. Obviously, we obtain that $K_{1,3}$ is an induced subgraph of $T$.

 If $n=4$, then $T\cong K_{1,3}$, and by Lemma 2.7, we have $\rho_{ag}(T)=2$.

 We now assume that $n\geq5$. Then $T$ contains $T_{1}$ or $T_{2}$ (as shown in Fig. 2) as its induced subgraph.

\begin{figure}[H]
\begin{center}
    \includegraphics[width=8cm]{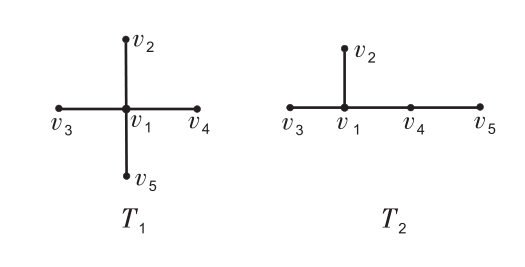} 
  \end{center}
  \caption*{\small \textbf{Fig. 2.} The trees $T_{1}$ and $T_{2}$.}
\end{figure}

\noindent {\textbf{Subcase 2.1}} If $T_{1}$ is an induced subgraph of $T$, then
$\displaystyle\rho_{ag}(T)\geq\rho(A(T))\geq\rho(A(T_{1}))$. Since

\begin{equation*}       
A(T_{1})=\left (                 
  \begin{array}{ccccc}   
    0 & 1 & 1 & 1 & 1\\  
    1 & 0 & 0 & 0 & 0\\  
    1 & 0 & 0 & 0 & 0\\
    1 & 0 & 0 & 0 & 0\\
    1 & 0 & 0 & 0 & 0\\
  \end{array}
\right).                 
\end{equation*}
A routine computation gives rise to $\rho(A(T_{1}))=2 $. Then $\displaystyle\rho_{ag}(T)\geq\rho(A(T_{1}))=2$.

\noindent {\textbf{Subcase 2.2}} If $T$ contains $T_{2}$ as its induced subgraph.

\noindent If $n=5$, then $T\cong T_{2}$ and

\begin{equation*}       
A_{ag}(T_{2})=\left (                 
  \begin{array}{ccccc}   
    0 & \frac{2\sqrt{3}}{3} & \frac{2\sqrt{3}}{3} & \frac{5\sqrt{6}}{12} & 0\\  
    \frac{2\sqrt{3}}{3} & 0 & 0 & 0 & 0\\  
    \frac{2\sqrt{3}}{3} & 0 & 0 & 0 & 0\\
    \frac{5\sqrt{6}}{12} & 0 & 0 & 0 & \frac{3\sqrt{2}}{4}\\
    0 & 0 & 0 & \frac{3\sqrt{2}}{4} & 0\\
  \end{array}
\right).                 
\end{equation*}

\noindent By computation, we obtain $\displaystyle\rho_{ag}(T)=\rho_{ag}(T_{2})\approx 2.0253>2$.

\noindent If $n\geq6$, then $T$ contains $T_{21}$ or $T_{22}$  (as shown in Fig. 3) as its induced subgraph.

\begin{figure}[H]
\begin{center}
    \includegraphics[width=12cm]{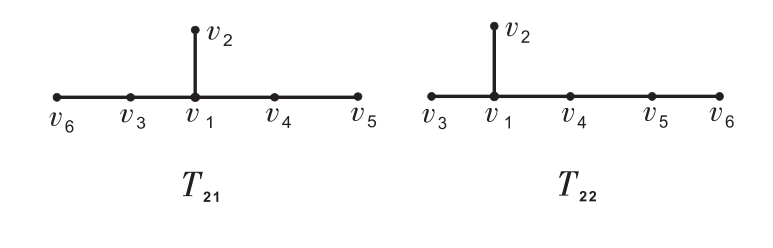} 
  \end{center}
  \caption*{\small \textbf{Fig. 3.} The trees $T_{21}$ and $T_{22}$.}
\end{figure}

\noindent \noindent {\textbf{Subcase 2.2.1}} $T_{21}$ is an induced subgraph of $T$. If $d_{2}=2$, then $T_{3}$, which is
depicted in Fig. 4, is an induced subgraph of $T$.

 \begin{figure}[H]
\begin{center}
    \includegraphics[width=5.3cm]{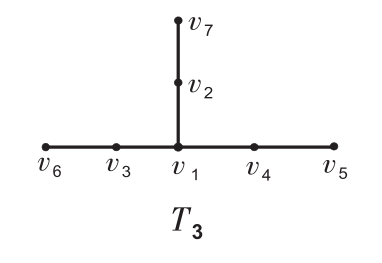} 
  \end{center}
  \caption*{\small \textbf{Fig. 4.} The tree $T_{3}$.}
\end{figure}

\noindent By directly calculation, we obtain
$\displaystyle\rho_{ag}(T)\geq\rho(A(T))\geq\rho(A(T_{3}))=2$.

If $d_{2}=1$, by $\displaystyle h(x,y)\geq 1=h(2,2)$ and Lemmas
2.2,2.3, we get $\rho_{ag}(T)\geq\rho(A_{ag}^{\prime}(T_{4}))$. The ``graph'' $T_{4}$  is shown in Fig. 5, where
$d_{5}=2$, $d_{6}=2$, $V(T_{4})=\{v_{1},v_{2},v_{3},v_{4},v_{5},v_{6}\}$ and $E(T_{4})=\{v_{1}v_{2}, v_{1}v_{3}, v_{1}v_{4},
v_{4}v_{5}, v_{3}v_{6}\}$. The matrix $A_{ag}^{\prime}(T_{4})=(h_{i,j})$, where $\displaystyle
h_{i,j}=\frac{d_{i}+d_{j}}{2\sqrt{d_{i}d_{j}}}$ if $v_{i}v_{j}\in E(T_{4})$ and 0 otherwise.

\begin{figure}[H]
\begin{center}
    \includegraphics[width=7cm]{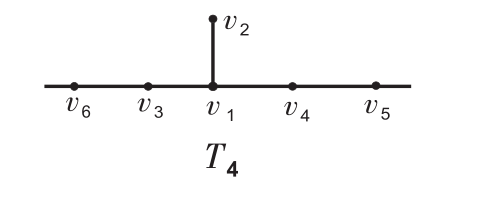} 
  \end{center}
  \caption*{\small \textbf{Fig. 5.} The tree $T_{4}$.}
\end{figure}

 \noindent  We have

\begin{equation*}       
A_{ag}^{\prime}(T_{4})=\left (                 
  \begin{array}{ccccccc}   
    0 & \frac{2\sqrt{3}}{3} & \frac{5\sqrt{6}}{12} & \frac{5\sqrt{6}}{12} & 0 & 0\\  
    \frac{2\sqrt{3}}{3} & 0 & 0 & 0 & 0 & 0\\  
    \frac{5\sqrt{6}}{12} & 0 & 0 & 0 & 0 & 1\\
    \frac{5\sqrt{6}}{12} & 0 & 0 & 0 & 1 & 0\\  
    0 & 0 & 0 & 1 & 0 & 0\\
    0 & 0 & 1 & 0 & 0 & 0\\  
  \end{array}
\right),                 
\end{equation*}

\noindent then
$\displaystyle\rho_{ag}(T)\geq\rho(A_{ag}^{\prime}(T_{4}))\thickapprox2.0226>2$.

\noindent {\textbf{Subcase 2.2.2}} $T_{22}$ is an induced subgraph of $T$.

If $d_{5}\geq3$, then $T_{5}$, which is depicted in Fig. 6, is an induced subgraph of $T$.

 \begin{figure}[H]
\begin{center}
    \includegraphics[width=7cm]{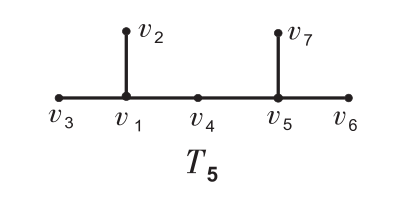} 
  \end{center}
  \caption*{\small \textbf{Fig. 6.} The tree $T_{5}$.}
\end{figure}

\noindent We acquire
$\displaystyle\rho_{ag}(T)\geq\rho(A(T))\geq\rho(A(T_{5}))=2$.

If $d_{5}=2$, thanks to $d_{1}\geq3$, then $d_{2},d_{3}=1$ or $2$, respectively.

If $d_{2}=2$ and $d_{3}=2$, then $T_{6}$(as shown in Fig. 7) is an induced subgraph of $T$.

 \begin{figure}[H]
\begin{center}
    \includegraphics[width=7cm]{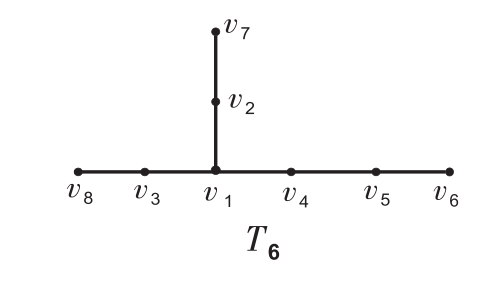} 
  \end{center}
  \caption*{\small \textbf{Fig. 7.} The tree $T_{6}$.}
\end{figure}

\noindent We obtain
$\displaystyle\rho_{ag}(T)\geq\rho(A(T))\geq\rho(A(T_{6}))\approx2.0285>2$.

If $d_{2}=1$ and $d_{3}=2$ or $d_{2}=2$ and $d_{3}=1$. Without loss of generality, we let $d_{2}=1$ and $d_{3}=2$. Due to $\displaystyle h(x,y)\geq 1=h(2,2)$ for $x,
y\geq1$ and Lemmas 2.2,2.3, we get $\rho_{ag}(T)\geq\rho(A_{ag}^{\prime}(T_{7}))$. The ``graph'' $T_{7}$  is shown in
Fig. 8, where $d_{6}=2$, $d_{7}=2$, $V(T_{7})=\{v_{1},v_{2},v_{3},v_{4},v_{5},v_{6},v_{7}\}$ and $E(T_{7})=\{v_{1}v_{2},
v_{1}v_{3}, v_{1}v_{4}, v_{4}v_{5}, v_{5}v_{6}, v_{3}v_{7}\}$. The matrix $A_{ag}^{\prime}(T_{7})=(h_{i,j})$, where $\displaystyle
h_{i,j}=\frac{d_{i}+d_{j}}{2\sqrt{d_{i}d_{j}}}$ if $v_{i}v_{j}\in E(T_{7})$ and 0 otherwise.

\begin{figure}[H]
\begin{center}
    \includegraphics[width=7cm]{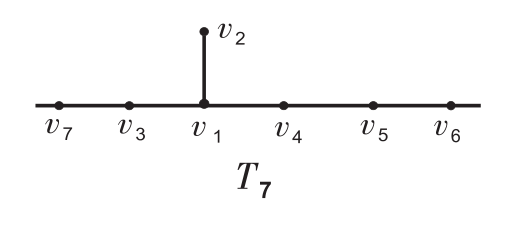} 
  \end{center}
  \caption*{\small \textbf{Fig. 8.} The tree $T_{7}$.}
\end{figure}

\noindent By simply calculation, we get
$\displaystyle\rho_{ag}(T)\geq\rho(A_{ag}^{\prime}(T_{7}))\thickapprox2.0523>2$.

 If $d_{2}=1$ and $d_{3}=1$, then
$\rho_{ag}(T)\geq\rho(A_{ag}^{\prime}(T_{8}))\thickapprox2.0457>2$, where the ``graph'' $T_{8}$  is shown in Fig. 9.

\begin{figure}[H]
\begin{center}
    \includegraphics[width=7cm]{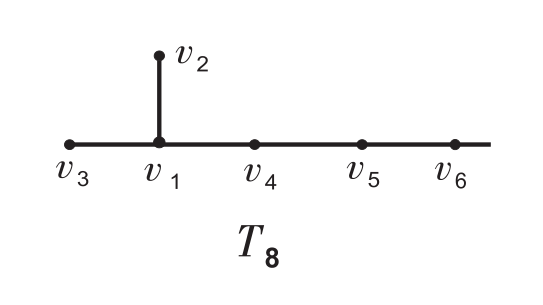} 
  \end{center}
  \caption*{\small \textbf{Fig. 9.} The tree $T_{8}$.}
\end{figure}

It remains to show that $\displaystyle2\cos\frac{\pi}{n+1}<\rho_{ag}(P_{n})<2$.
 The arithmetic-geometric matrix of $P_{n}$ can be written as

\begin{equation*}       
A_{ag}(P_{n})=                
  \begin{pmatrix}
    0 & \frac{3\sqrt{2}}{4} & 0 & 0 & \ldots & 0 & 0 & 0 & 0\\  
    \frac{3\sqrt{2}}{4} & 0 & 1 & 0 & \ldots & 0 & 0 & 0 & 0\\  
    0 & 1 & 0 & 1 & \ldots & 0 & 0 & 0 & 0\\
    0 & 0 & 1 & 0 & \ldots & 0 & 0 & 0 & 0\\  
    \vdots & \vdots & \vdots & \vdots & \ddots & \vdots & \vdots & \vdots & \vdots\\
    0 & 0 & 0 & 0 & \ldots & 0 & 1 & 0 & 0\\  
    0 & 0 & 0 & 0 & \ldots & 1 & 0 & 1 & 0\\
    0 & 0 & 0 & 0 & \ldots & 0 & 1 & 0 & \frac{3\sqrt{2}}{4}\\
    0 & 0 & 0 & 0 & \ldots & 0 & 0 & \frac{3\sqrt{2}}{4} & 0\\

  \end{pmatrix}.                
\end{equation*}

\noindent
\begin{equation*}       
f(\rho,n)=det(\rho E_{n}-A_{ag}(P_{n}))=                
  \begin{vmatrix}
   \rho & -\frac{3\sqrt{2}}{4} & 0 & 0 & \ldots & 0 & 0 & 0 & 0\\  
    -\frac{3\sqrt{2}}{4} &\rho & -1 & 0 & \ldots & 0 & 0 & 0 & 0\\  
    0 & -1 &\rho & -1 & \ldots & 0 & 0 & 0 & 0\\
    0 & 0 & -1 &\rho & \ldots & 0 & 0 & 0 & 0\\  
    \vdots & \vdots & \vdots & \vdots & \ddots & \vdots & \vdots & \vdots & \vdots\\
    0 & 0 & 0 & 0 & \ldots &\rho & -1 & 0 & 0\\  
    0 & 0 & 0 & 0 & \ldots & -1 &\rho & -1 & 0\\
    0 & 0 & 0 & 0 & \ldots & 0 & -1 &\rho & -\frac{3\sqrt{2}}{4}\\
    0 & 0 & 0 & 0 & \ldots & 0 & 0 & -\frac{3\sqrt{2}}{4} &\rho\\
  \end{vmatrix}.                
\end{equation*}

\noindent If $\displaystyle\rho=2$, by directly calculation, we obtain
$\displaystyle f(2,n)=\frac{49n}{64}+\frac{77}{64}>0.$

\noindent If $\displaystyle\rho>2$, then by Lemma 2.6, we obtain
\begin{eqnarray*}
f(\rho,n)&=&\rho^{2}\Phi(P_{n-2};\rho)-\frac{9}{4}\rho\Phi(P_{n-3};\rho)+\Big(\frac{9}{8}\Big)^{2}\Phi(P_{n-4};\rho)\\
&=&\frac{1}{2^{n-3}\sqrt{\rho^{2}-4}}\Big\{\Big(\frac{\rho}{2}\Big)^{2}\Big[\Big(\rho+\sqrt{\rho^{2}-4}\Big)^{n-1}-\Big(\rho-\sqrt{\rho^{2}-4}\Big)^{n-1}\Big]-\frac{9\rho}{8}\Big[\Big(\rho+\sqrt{\rho^{2}-4}\Big)^{n-2}\\
&-&\Big(\rho-\sqrt{\rho^{2}-4}\Big)^{n-2}\Big]+\Big(\frac{9}{8}\Big)^{2}\Big[\Big(\rho+\sqrt{\rho^{2}-4}\Big)^{n-3}-\Big(\rho-\sqrt{\rho^{2}-4}\Big)^{n-3}\Big]\Big\}.
\end{eqnarray*}

 Now it suffices to prove that
\begin{eqnarray*}
 g(\rho,n)&=&\Big(\frac{\rho}{2}\Big)^{2}\Big[\Big(\rho+\sqrt{\rho^{2}-4}\Big)^{n-1}-\Big(\rho-\sqrt{\rho^{2}-4}\Big)^{n-1}\Big]
 -\frac{9\rho}{8}\Big[\Big(\rho+\sqrt{\rho^{2}-4}\Big)^{n-2}\\
&-&\Big(\rho-\sqrt{\rho^{2}-4}\Big)^{n-2}\Big]+\Big(\frac{9}{8}\Big)^{2}\Big[\Big(\rho+\sqrt{\rho^{2}-4}\Big)^{n-3}-\Big(\rho-\sqrt{\rho^{2}-4}\Big)^{n-3}\Big]>0
\end{eqnarray*}for $\rho>2$.

  If  $\displaystyle
 \frac{\rho(\rho-\sqrt{\rho^{2}-4})}{2}-\frac{9}{8}>0$, then

 \begin{center}
 $\displaystyle\Big[\frac{\rho(\rho+\sqrt{\rho^{2}-4})}{2}-\frac{9}{8}\Big]^{2}>\Big[\frac{\rho(\rho-\sqrt{\rho^{2}-4})}{2}-\frac{9}{8}\Big]^{2}$.
\end{center}

 \noindent  If $\displaystyle \frac{\rho(\rho-\sqrt{\rho^{2}-4})}{2}-\frac{9}{8}\leq0$,
then

 \begin{center}
 $\displaystyle\frac{\rho(\rho+\sqrt{\rho^{2}-4})}{2}-\frac{9}{8}+\frac{\rho(\rho-\sqrt{\rho^{2}-4})}{2}-\frac{9}{8}=\rho^{2}-\frac{9}{4}>0$.
\end{center}

\noindent Thus we also have $\displaystyle\Big[\frac{\rho(\rho+\sqrt{\rho^{2}-4})}{2}-\frac{9}{8}\Big]^{2}>\Big[\frac{\rho(\rho-\sqrt{\rho^{2}-4})}{2}-\frac{9}{8}\Big]^{2}$.

 \noindent  Then for $\displaystyle\rho>2$, we obtain

\begin{center}
\noindent $\displaystyle(\rho+\sqrt{\rho^{2}-4})^{n-3}\Big[\frac{\rho(\rho+\sqrt{\rho^{2}-4})}{2}-\frac{9}{8}\Big]^{2}$
\end{center}

$ \displaystyle\ \ \ \ \ \ \ \ \ \ \ \ \ \ \ \ \ \ \ \ \ \ \ \ \ \ \ \ \ \
>(\rho-\sqrt{\rho^{2}-4})^{n-3}\Big[\frac{\rho(\rho+\sqrt{\rho^{2}-4})}{2}-\frac{9}{8}\Big]^{2}$

$\displaystyle\ \ \ \ \ \ \ \ \ \ \ \ \ \ \ \ \ \ \ \ \ \ \ \ \ \ \ \ \ \
>(\rho-\sqrt{\rho^{2}-4})^{n-3}\Big[\frac{\rho(\rho-\sqrt{\rho^{2}-4})}{2}-\frac{9}{8}\Big]^{2}$.

\noindent It follows that

\begin{center}
$\displaystyle(\rho+\sqrt{\rho^{2}-4})^{n-3}\Big[\Big(\frac{\rho}{2}\Big)^{2}(\rho+\sqrt{\rho^{2}-4})^{2}-\frac{9\rho}{8}(\rho+\sqrt{\rho^{2}-4})+\Big(\frac{9}{8}\Big)^{2}\Big]$
\end{center}

$\displaystyle\ \ \ \ \ \ \
>(\rho-\sqrt{\rho^{2}-4})^{n-3}\Big[\Big(\frac{\rho}{2}\Big)^{2}(\rho-\sqrt{\rho^{2}-4})^{2}-\frac{9\rho}{8}(\rho-\sqrt{\rho^{2}-4})+\Big(\frac{9}{8}\Big)^{2}\Big]$.

\noindent That is $\displaystyle g(\rho,n)>0$.

 Because for $P_{n}$, $\displaystyle A_{ag}(P_{n})> A(P_{n})$, then we have $\displaystyle \rho_{ag}(P_{n})>\rho(A(P_{n}))=2\cos\frac{\pi}{n+1}$.

 \hfill$\Box$

 \section{The $AG$ spectral radius of unicyclic graphs}

\noindent {\textbf{Theorem 4.1.}}
Let $G$ be an unicyclic of order $n\geq3$, then

\begin{center}
$\displaystyle 2=\rho_{ag}(C_{n})\leq\rho_{ag}(G)\leq\rho_{ag}(S_{n}+e)$
\end{center}

\noindent with equality if and only if $T\cong C_{n}$ for the lower bound, and if and only if $T\cong S_{n}+e$ for the upper bound and
$\displaystyle\rho_{ag}(S_{n}+e)<\frac{n}{2}$ for $n\geq7$.

\noindent {\textbf{Proof.}} If $n=3$, it is manifest. For any unicyclic graph $G$ of order $n\geq4$, from Lemmas 2.2, 2.5 and
$\displaystyle h(x,y)\geq h(2,2)=1$ for $x,y\geq1$, we
have $\displaystyle\rho_{ag}(G)\geq\rho(A(G)) \geq\rho(A(C_{n}))=\rho_{ag}(C_{n})=2.$ This leads to the lower bound.

Our task now is to consider the upper bound. The approximate values of $\rho_{ag}$ of unicyclic graphs with order $4\leq n\leq7$ are given in Tables 1-4 which are shown in the  Appendix 1. It is obvious that
$S_{n}+e$ has the largest $AG$ spectral radius and $\displaystyle\rho_{ag}(S_{7}+e)<\frac{n}{2}$.

For $n\geq8$, according to Lemma 2.1, if we find a vector $\textbf{x}>\textbf{0}$ such that $\displaystyle A_{ag}(G)\textbf{x}\leq
\frac{n-1}{2}\textbf{x}$, then
$\displaystyle\rho_{ag}(G)\leq\frac{n-1}{2}$.
Let $x_{i}=\sqrt{d_{i}}$.

Suppose that any pair of adjacent vertices $v_{i},v_{j}\in G$, we have $d_{i}+d_{j}\leq n-1$, then for any vertex $v_{i}$,

$\displaystyle (A_{ag}(G)\textbf{x})_{i}=\sum_{j=1,i\sim j}^{n}\frac{d_{i}+d_{j}}{2\sqrt{d_{i}d_{j}}}\cdot \sqrt{d_{j}}
 =\sum_{j=1,i\sim j}^{n}\frac{d_{i}+d_{j}}{2\sqrt{d_{i}}}\leq \frac{n-1}{2\sqrt{d_{i}}}\cdot d_{i}=\frac{n-1}{2}\cdot\sqrt{d_{i}}$.

\noindent  It follows that $\displaystyle\rho_{ag}(G)\leq\frac{n-1}{2}$.

Next we assume that there exist a pair of adjacent vertices $v_{i},v_{j}$ such that $d_{i}+d_{j}\geq n$.
By the proof of Conjecture 1 in \cite{13}, we know if $G$ exists two adjacent vertices $v_{i},v_{j}$ such that $d_{i}+d_{j}\geq n$,
then $G\in\{\mathcal{G}_{1},\mathcal{G}_{2},\mathcal{G}_{3},\mathcal{G}_{4},\mathcal{G}_{5}\}$(as shown in Fig. 10).

\begin{figure}[H]
\begin{center}
    \includegraphics[width=12cm]{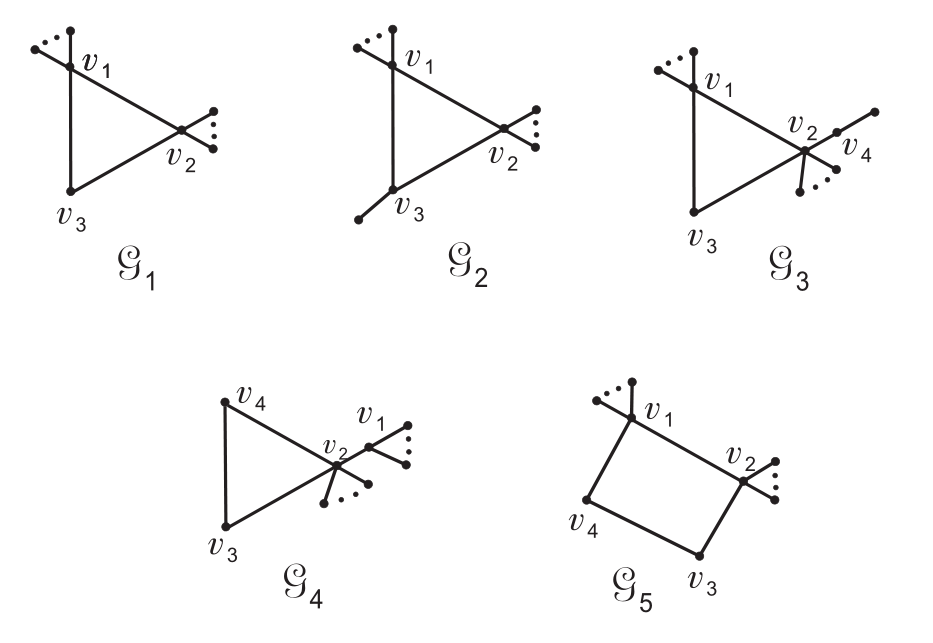} 
  \end{center}
  \caption*{\small \textbf{Fig. 10.} The graphs  $\mathcal{G}_{1},\mathcal{G}_{2},\mathcal{G}_{3},\mathcal{G}_{4},\mathcal{G}_{5}$.}
\end{figure}

The graphs ${G}_{1},{G}_{2},{G}_{3},{G}_{4}$ shown in Fig. 11 will be used in the following.

\begin{figure}[H]
\begin{center}
    \includegraphics[width=13cm]{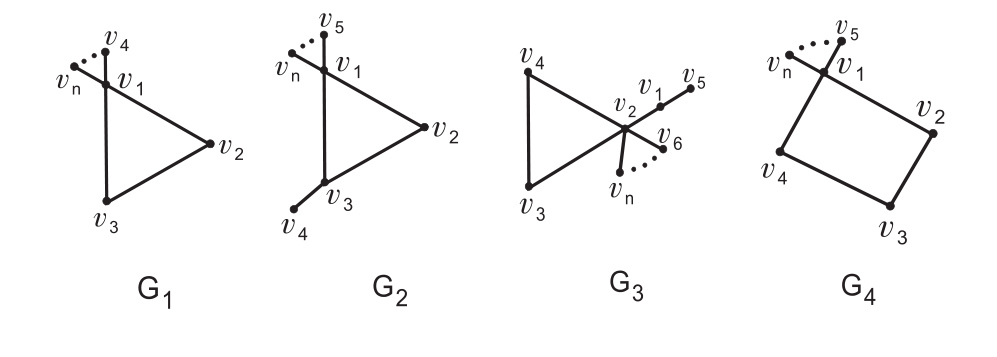} 
  \end{center}
  \caption*{\small \textbf{Fig. 11.} The graphs  ${G}_{1},{G}_{2},{G}_{3},{G}_{4}$.}
\end{figure}

Let $G\in\mathcal{G}_{1}$ and $G\ncong G_{1},G_{2}$, then $d_{1}+d_{2}=n+1$ and $4\leq
d_{1},d_{2}\leq
n-3$.

\noindent For the pendent vertex $v$,

\begin{center}
$\displaystyle
(A_{ag}(G)\textbf{x})_{v}=\frac{1+d_{i}}{2\sqrt{d_{i}}}\cdot\sqrt{d_{i}}=\frac{1+d_{i}}{2}<\frac{n-1}{2}$ for $i=1,2$.
\end{center}

\noindent For the vertex $v_{1}$,
\begin{eqnarray*}
(A_{ag}(G)\textbf{x})_{v_{1}}&=&(d_{1}-2)\cdot\frac{1+d_{1}}{2\sqrt{d_{1}}}+\frac{d_{1}+d_{2}}{2\sqrt{d_{1}d_{2}}}\cdot\sqrt{d_{2}}+\frac{d_{1}+2}{2\sqrt{2d_{1}}}\cdot\sqrt{2}
\\
&=&\frac{(d_{1}-2)(1+d_{1})}{2\sqrt{d_{1}}}+\frac{d_{1}+d_{2}}{2\sqrt{d_{1}}}+\frac{d_{1}+2}{2\sqrt{d_{1}}}\\
&=&\frac{(d_{1}-2)(1+d_{1})+2d_{1}+d_{2}+2}{2\sqrt{d_{1}}}=\frac{d_{1}^{2}+n+1}{2\sqrt{d_{1}}}\\
&\leq& \frac{(n-1)d_{1}}{2\sqrt{d_{1}}}=\frac{n-1}{2}\cdot\sqrt{d_{1}}.
\end{eqnarray*}

\noindent Similarly, for the vertex $v_{2}$, we have
\begin{center}
$\displaystyle (A_{ag}(G)\textbf{x})_{v_{2}}\leq \frac{(n-1)}{2}\cdot\sqrt{d_{2}}$.
\end{center}

\noindent For the vertex $v_{3}$,

$\displaystyle
(A_{ag}(G)\textbf{x})_{v_{3}}=\frac{2+d_{1}}{2\sqrt{2d_{1}}}\cdot\sqrt{d_{1}}+\frac{2+d_{2}}{2\sqrt{2d_{2}}}\cdot\sqrt{d_{2}}=\frac{4+d_{1}+d_{2}}{2\sqrt{2}}=\frac{n+5}{2\sqrt{2}}\leq
\frac{n-1}{2}\cdot\sqrt{2}$.

Thus $\displaystyle (A_{ag}(G)\textbf{x})_{i}\leq \frac{n-1}{2}x_{i}$ holds for all vertices in $G\in\mathcal{G}_{1}$ and $G\ncong
G_{1},G_{2}$.

The proof of $G\in\mathcal{G}_{2}$ and $G\ncong
G_{2}$, $G\in\mathcal{G}_{3}$ and $G\ncong
G_{3}$, $G\in\mathcal{G}_{4}$ and $G\ncong
G_{3}$ or $G\in\mathcal{G}_{5}$ and $G\ncong
G_{5}$ is quite similar to
that given for the case of $G\in\mathcal{G}_{1}$ and $G\ncong
G_{1},G_{2}$ and given in the Appendix 1.

 Then it remains to prove that
$\displaystyle\rho_{ag}(G_{2}),\rho_{ag}(G_{3}),\rho_{ag}(G_{4})<\frac{n-1}{2}<\rho_{ag}(G_{1})<\frac{n}{2}$.

\begin{equation*}       
  A_{ag}(G_{1})=\begin{pmatrix}
    0 & \frac{n+1}{2\sqrt{2(n-1)}} & \frac{n+1}{2\sqrt{2(n-1)}} & \frac{n}{2\sqrt{n-1}} & \ldots & \frac{n}{2\sqrt{n-1}} &
    \frac{n}{2\sqrt{n-1}}\\  
    \frac{n+1}{2\sqrt{2(n-1)}} & 0 & 1 & 0 & \ldots & 0 & 0\\  
    \frac{n+1}{2\sqrt{2(n-1)}} & 1 & 0 & 0 & \ldots & 0 & 0\\
    \frac{n}{2\sqrt{n-1}} & 0 & 0 & 0 & \ldots & 0 & 0\\  
    \vdots & \vdots & \vdots & \vdots & \ddots & \vdots & \vdots\\
    \frac{n}{2\sqrt{n-1}} & 0 & 0 & 0 & \ldots & 0 & 0\\
    \frac{n}{2\sqrt{n-1}} & 0 & 0 & 0 & \ldots & 0 & 0\\

  \end{pmatrix}.                
\end{equation*}

\noindent Its arithmetic-geometric characteristic polynomial
$g_{1}(\rho,n)=det(\rho E_{n}-A_{ag}(G_{1}))$

\begin{equation*}       
=                
  \begin{vmatrix}
   \rho & -\frac{n+1}{2\sqrt{2(n-1)}} & -\frac{n+1}{2\sqrt{2(n-1)}} & -\frac{n}{2\sqrt{n-1}} & \ldots & -\frac{n}{2\sqrt{n-1}} &
   -\frac{n}{2\sqrt{n-1}}\\  
    -\frac{n+1}{2\sqrt{2(n-1)}} &\rho & -1 & 0 & \ldots & 0 & 0\\  
    -\frac{n+1}{2\sqrt{2(n-1)}} & -1 &\rho & 0 & \ldots & 0 & 0\\
    -\frac{n}{2\sqrt{n-1}} & 0 & 0 &\rho & \ldots & 0 & 0\\  
    \vdots & \vdots & \vdots & \vdots & \ddots & \vdots & \vdots\\
    -\frac{n}{2\sqrt{n-1}} & 0 & 0 & 0 & \ldots &\rho & 0\\
    -\frac{n}{2\sqrt{n-1}} & 0 & 0 & 0 & \ldots & 0 &\rho\\

  \end{vmatrix}              
\end{equation*}

$\displaystyle \ \ \ \ \ \ \ \
=\rho^{n-4}(\rho+1)\Big[\rho^{3}-\rho^{2}-\frac{n^{3}-2n^{2}+2n+1}{4(n-1)}\rho+\frac{(n-3)n^{2}}{4(n-1)}\Big]$.

Let $\displaystyle t_{1}(\rho,n)=\rho^{3}-\rho^{2}-\frac{n^{3}-2n^{2}+2n+1}{4(n-1)}\rho+\frac{(n-3)n^{2}}{4(n-1)}$, denoted by
$\displaystyle\rho_{1}\geq\rho_{2}\geq\rho_{3}$ the three roots of $\displaystyle t_{1}(\rho,n)=0$. Obviously,
$\rho_{1}=\rho_{ag}$, i.e., the desired result is equivalent to $\displaystyle \frac{n-1}{2}<\rho_{1}<\frac{n}{2}$.
From $\displaystyle t_{1}(0,n)=\frac{(n-3)n^{2}}{4(n-1)}>0$,

$\displaystyle t_{1}\Big(\frac{n-1}{2},n\Big)=\frac{(n-1)^{4}-2(n-1)^{3}-(n-1)(n^{3}-2n^{2}+2n+1)+2n^{2}(n-3)}{8(n-1)}$

$\displaystyle \ \ \ \ \ \ \ \ \ \ \ \ \ \ \ \ \ \ \ \ \ =\frac{-n^{3}+2n^{2}-9n+4}{8(n-1)}<0$ for $n\geq8$.

\noindent and

$\displaystyle t_{1}\Big(\frac{n}{2},n\Big)=\frac{(n-1)n^{3}-2n^{2}(n-1)-n(n^{3}-2n^{2}+2n+1)+2n^{2}(n-3)}{8(n-1)}$

$\displaystyle \ \ \ \ \ \ \ \ \ \ \ \ \  \ =\frac{n^{3}-6n^{2}-n}{8(n-1)}>0$ for $n\geq8$.

\noindent  It follows that $\displaystyle \frac{n-1}{2}<\rho_{1}<\frac{n}{2}$, i.e.,
$\displaystyle \frac{n-1}{2}<\rho_{ag}(G_{1})<\frac{n}{2}$.

Similarly we obtain the arithmetic-geometric characteristic polynomial $g_{i}(\rho,n)$ of $G_{i} (i=2,3,4)$ as below
\end{spacing}
\begin{spacing}{1.6}
\noindent $(1)$  $\displaystyle
g_{2}(\rho,n)=\frac{\rho^{n-4}}{4(n-2)}\Big[4(n-2)\rho^{4}-\frac{6n^{3}-31n^{2}+115n-136}{6}\rho^{2}-\frac{5n^{2}+5n}{6}\rho+\frac{2n^{2}}{3}$

$\displaystyle \ \ \ \ \ \ \ \ \ \ \ \ \ \ \ \ \ +\frac{19(n-4)(n-1)^{2}}{8}\Big]$;
\end{spacing}
\begin{spacing}{1.6}

\end{spacing}
\begin{spacing}{1.5}
\noindent $(2)$ $\displaystyle
 g_{3}(\rho,n)=\frac{\rho^{n-6}}{8(n-2)}\Big\{8(n-2)\rho^{6}-(2n^{3}-11n^{2}+39n-44)\rho^{4}-2n^{2}\rho^{3}+\Big[\frac{13n^{2}}{4}+9(n-2)+$

 $\displaystyle \ \ \ \ \ \ \ \ \ \ \ \ \ \ \ \ \
 \frac{17}{4}(n-5)(n-1)^{2}\Big]\rho^{2}+\frac{9}{4}n^{2}\rho-\frac{9}{4}(n-5)(n-1)^{2}\Big\}$;
\end{spacing}
\begin{spacing}{1.6}

\end{spacing}
\begin{spacing}{1.5}

\noindent $(3)$ $\displaystyle
g_{4}(\rho,n)=\frac{\rho^{n-4}}{8(n-2)}[8(n-2)\rho^{4}-(2n^{3}-10n^{2}+34n-40)\rho^{2}+4(n-4)(n-1)^{2}]$.
\end{spacing}
\begin{spacing}{2}

For $\displaystyle \rho\geq\frac{n-1}{2}$ and $n\geq14$, we acquire

\noindent (1) $\displaystyle (4n-9)\rho^{4}>\frac{6n^{3}-31n^{2}+115n-136}{6}\rho^{2}$ and
$\displaystyle \rho^{4}>\frac{5n^{2}+5n}{6}\rho$, then $\displaystyle g_{2}(\rho,n)>0$;

\noindent (2)  $\displaystyle (8n-18)\rho^{6}>(2n^{3}-11n^{2}+39n-44)\rho^{4}$, $\displaystyle
9(n-2)\rho^{2}+\frac{9}{4}n^{2}\rho>\frac{9}{4}(n-5)(n-1)^{2}$ and

$\displaystyle
\rho^{6}>n^{2}\rho^{3}$, then $\displaystyle g_{3}(\rho,n)>0$;

\noindent (3)  $\displaystyle 8(n-2)\rho^{4}>(2n^{3}-10n^{2}+34n-40)\rho^{2}$, then $\displaystyle
g_{4}(\rho,n)>0$.

For $8\leq n\leq13$, we can easily vertify that $\displaystyle g_{2}(\rho,n), g_{3}(\rho,n)$ and $\displaystyle g_{4}(\rho,n)>0$ for
$\displaystyle \rho\geq\frac{n-1}{2}$.

  \hfill$\Box$

 \end{spacing}

\section*{Acknowledgments}

\ \ \ \ \ This work is supported by NSFC (No. 11671336) and the Fundamental Research Funds for the Central Universities (No. 20720190062).

\bibliographystyle{abbrv}
\bibliography{References}

\newpage
 \begin{center}
$\displaystyle \bf Appendix$ $\displaystyle 1$
\end{center}
\begin{table}[h]
\renewcommand\arraystretch{2.5}
  \centering
  \begin{tabular}{ | c | l | l |l |}
    \hline
    The unicyclic graph G & $\displaystyle\rho_{ag}(G)$ & The unicyclic graph G & $\displaystyle\rho_{ag}(G)$\\
    \hline
    \begin{minipage}[b]{0.3\columnwidth}
		\centering
		\raisebox{-.30\height}{\includegraphics[width=1.8cm]{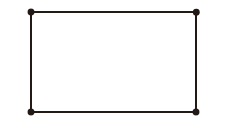}}
	\end{minipage}
    & \ \ \ \ \ 2
    & \begin{minipage}[b]{0.3\columnwidth}
		\centering
		\raisebox{-.30\height}{\includegraphics[width=2cm]{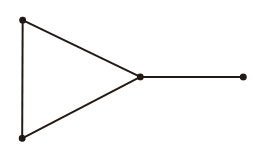}}
	\end{minipage}
    & 2.2536
    \\ [3pt]\hline
 \end{tabular}
  \caption{The approximate value of $\rho_{ag}$ of the unicyclic graph with order $n=4$}
\end{table}

   \begin{table}[h]
\renewcommand\arraystretch{2.5}
  \centering
  \begin{tabular}{ | c | l | l |l | }
    \hline
The unicyclic graph G & $\displaystyle\rho_{ag}(G)$ & The unicyclic graph G & $\displaystyle\rho_{ag}(G)$\\
    \hline
    \begin{minipage}[b]{0.3\columnwidth}
		\centering
		\raisebox{-.35\height}{\includegraphics[width=1.4cm]{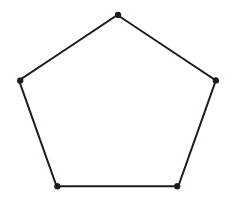}}
	\end{minipage}
    & \ \ \ \ \ 2
    & \begin{minipage}[b]{0.3\columnwidth}
		\centering
		\raisebox{-.35\height}{\includegraphics[width=1.75cm]{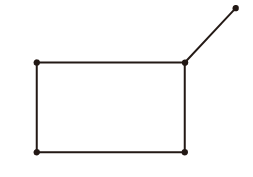}}
	\end{minipage}
    &  2.2066
    \\[3pt]
    \hline
   \begin{minipage}[b]{0.3\columnwidth}
		\centering
		\raisebox{-.30\height}{\includegraphics[width=2.3cm]{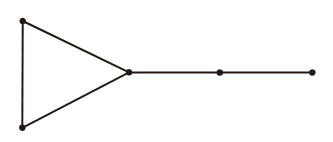}}
	\end{minipage}
    & 2.2543
    & \begin{minipage}[b]{0.3\columnwidth}
		\centering
		\raisebox{-.35\height}{\includegraphics[width=1.5cm]{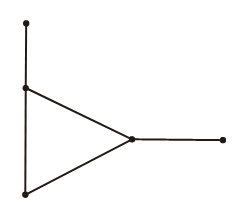}}
	\end{minipage}
    &  2.4149
    \\[3pt]
    \hline
  \begin{minipage}[b]{0.3\columnwidth}
		\centering
		\raisebox{-.30\height}{\includegraphics[width=1.7cm]{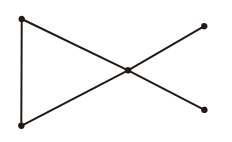}}
	\end{minipage}
    & 2.6035
    &
    &
    \\ \hline
\end{tabular}
  \caption{The approximate value of $\rho_{ag}$ of the unicyclic graph with order $n=5$}
\end{table}

\begin{table}[h]
\renewcommand\arraystretch{2.5}
  \centering
  \begin{tabular}{ | c | l | l |l |}
    \hline
   The unicyclic graph G & $\displaystyle\rho_{ag}(G)$ & The unicyclic graph G & $\displaystyle\rho_{ag}(G)$\\
    \hline
    \begin{minipage}[b]{0.3\columnwidth}
		\centering
		\raisebox{-.25\height}{\includegraphics[width=1cm]{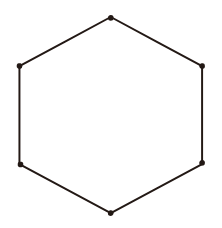}}
	\end{minipage}
    & \ \ \ \ 2
    & \begin{minipage}[b]{0.3\columnwidth}
		\centering
		\raisebox{-.28\height}{\includegraphics[width=1cm]{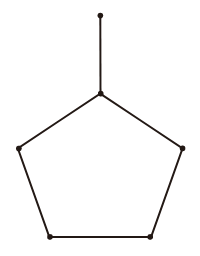}}
	\end{minipage}
    &  2.1785
    \\[3pt]
    \hline
    \begin{minipage}[b]{0.3\columnwidth}
		\centering
		\raisebox{-.28\height}{\includegraphics[width=1.5cm]{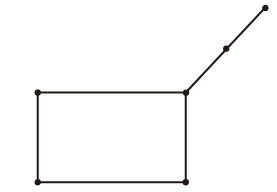}}
	\end{minipage}
    & 2.2096
    & \begin{minipage}[b]{0.3\columnwidth}
		\centering
		\raisebox{-.28\height}{\includegraphics[width=2cm]{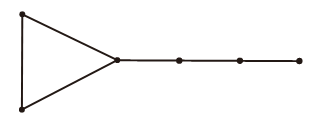}}
	\end{minipage}
    &  2.2632
    \\[3pt]
    \hline
 \begin{minipage}[b]{0.3\columnwidth}
		\centering
		\raisebox{-.33\height}{\includegraphics[width=1.5cm]{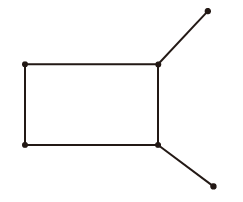}}
	\end{minipage}
    & 2.3439
    & \begin{minipage}[b]{0.3\columnwidth}
		\centering
		\raisebox{-.30\height}{\includegraphics[width=1.45cm]{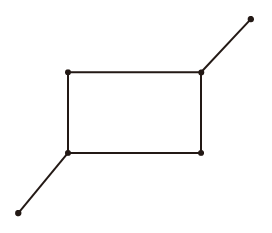}}
	\end{minipage}
    &  2.3452
    \\[3pt]
    \hline
     \begin{minipage}[b]{0.3\columnwidth}
		\centering
		\raisebox{-.25\height}{\includegraphics[width=1.6cm]{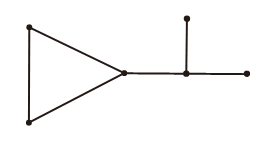}}
	\end{minipage}
    & 2.3551
    & \begin{minipage}[b]{0.3\columnwidth}
		\centering
		\raisebox{-.25\height}{\includegraphics[width=1.5cm]{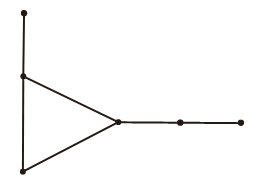}}
	\end{minipage}
    &  2.4095
    \\[3pt]
    \hline
     \begin{minipage}[b]{0.3\columnwidth}
		\centering
		\raisebox{-.35\height}{\includegraphics[width=1.2cm]{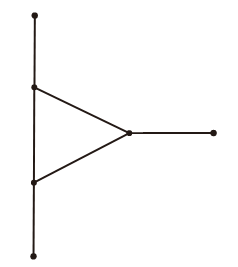}}
	\end{minipage}
    & 2.5275
    & \begin{minipage}[b]{0.3\columnwidth}
		\centering
		\raisebox{-.25\height}{\includegraphics[width=1.5cm]{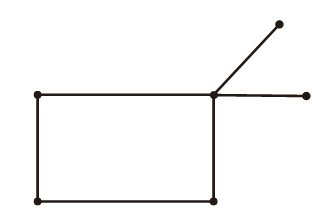}}
	\end{minipage}
    &  2.5295
    \\[3pt]
    \hline
     \begin{minipage}[b]{0.3\columnwidth}
		\centering
		\raisebox{-.25\height}{\includegraphics[width=1.6cm]{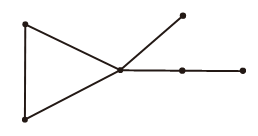}}
	\end{minipage}
    & 2.5695
    & \begin{minipage}[b]{0.3\columnwidth}
		\centering
		\raisebox{-.33\height}{\includegraphics[width=1.3cm]{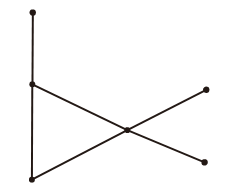}}
	\end{minipage}
    &  2.6879
    \\[3pt]
    \hline
     \begin{minipage}[b]{0.3\columnwidth}
		\centering
		\raisebox{-.30\height}{\includegraphics[width=1.3cm]{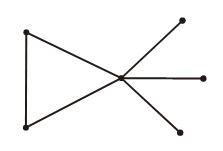}}
	\end{minipage}
    &  3.0113
    &
    &
     \\[3pt]
 \hline
\end{tabular}
  \caption{The approximate value of $\rho_{ag}$ of the unicyclic graph with order $n=6$}
\end{table}

\begin{table}[h]
\scalebox{0.1}{}
\renewcommand\arraystretch{2.2}
  \centering
  \begin{tabular}{ | c | l | l |l |}
    \hline
   The unicyclic graph G & $\displaystyle\rho_{ag}(G)$ & The unicyclic graph G & $\displaystyle\rho_{ag}(G)$\\
    \hline
     \begin{minipage}[b]{0.3\columnwidth}
		\centering
		\raisebox{-.25\height}{\includegraphics[width=0.9cm]{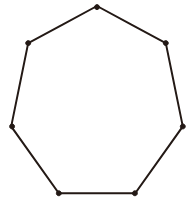}}
	\end{minipage}
    & \ \ \ \ 2
    & \begin{minipage}[b]{0.3\columnwidth}
		\centering
		\raisebox{-.25\height}{\includegraphics[width=0.7cm]{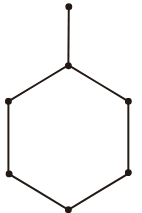}}
	\end{minipage}
    &  2.1602
    \\[2.5pt]
    \hline
     \begin{minipage}[b]{0.3\columnwidth}
		\centering
		\raisebox{-.25\height}{\includegraphics[width=1.5cm]{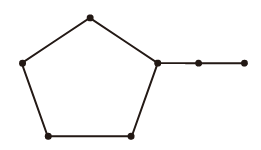}}
	\end{minipage}
    & 2.1827
    & \begin{minipage}[b]{0.3\columnwidth}
		\centering
		\raisebox{-.25\height}{\includegraphics[width=1.5cm]{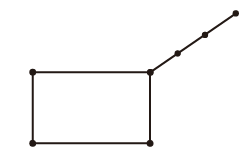}}
	\end{minipage}
    &  2.2188
    \\[2.5pt]
    \hline
    \begin{minipage}[b]{0.3\columnwidth}
		\centering
		\raisebox{-.25\height}{\includegraphics[width=1.6cm]{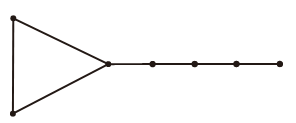}}
	\end{minipage}
    & 2.2661
    & \begin{minipage}[b]{0.3\columnwidth}
		\centering
		\raisebox{-.25\height}{\includegraphics[width=1.1cm]{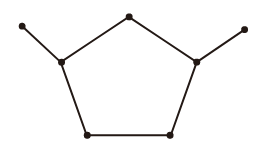}}
	\end{minipage}
    &  2.2942
    \\[2.5pt]
    \hline
     \begin{minipage}[b]{0.3\columnwidth}
		\centering
		\raisebox{-.30\height}{\includegraphics[width=1.2cm]{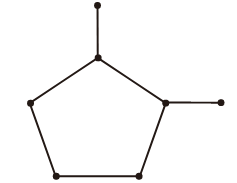}}
	\end{minipage}
    & 2.3041
    & \begin{minipage}[b]{0.3\columnwidth}
		\centering
		\raisebox{-.25\height}{\includegraphics[width=1.1cm]{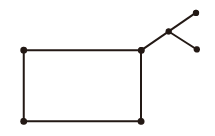}}
	\end{minipage}
    &  2.3094
    \\[2.5pt]
    \hline
     \begin{minipage}[b]{0.3\columnwidth}
		\centering
		\raisebox{-.25\height}{\includegraphics[width=1.7cm]{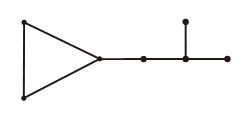}}
	\end{minipage}
    & 2.316
    & \begin{minipage}[b]{0.3\columnwidth}
		\centering
		\raisebox{-.25\height}{\includegraphics[width=1.1cm]{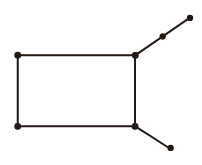}}
	\end{minipage}
    &  2.3413
    \\[2.5pt]
    \hline
    \begin{minipage}[b]{0.3\columnwidth}
		\centering
		\raisebox{-.25\height}{\includegraphics[width=1.4cm]{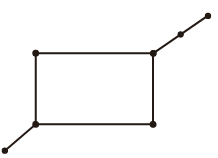}}
	\end{minipage}
    & 2.3428
    & \begin{minipage}[b]{0.3\columnwidth}
		\centering
		\raisebox{-.25\height}{\includegraphics[width=1.7cm]{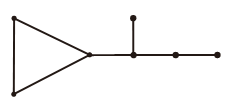}}
	\end{minipage}
    &  2.3528
    \\[2.5pt]
    \hline
    \begin{minipage}[b]{0.3\columnwidth}
		\centering
		\raisebox{-.30\height}{\includegraphics[width=1.2cm]{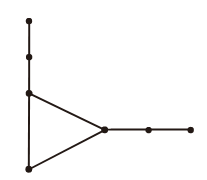}}
	\end{minipage}
    & 2.4044
    & \begin{minipage}[b]{0.3\columnwidth}
		\centering
		\raisebox{-.25\height}{\includegraphics[width=1.6cm]{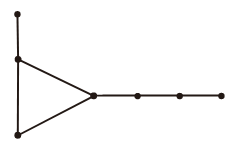}}
	\end{minipage}
    &  2.4135
    \\[2.5pt]
   \hline
   \begin{minipage}[b]{0.3\columnwidth}
		\centering
		\raisebox{-.25\height}{\includegraphics[width=1.3cm]{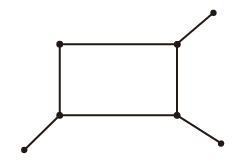}}
	\end{minipage}
    & 2.4473
    & \begin{minipage}[b]{0.3\columnwidth}
		\centering
		\raisebox{-.25\height}{\includegraphics[width=1.3cm]{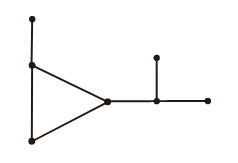}}
	\end{minipage}
    &  2.4732
    \\[2.5pt]
    \hline
    \begin{minipage}[b]{0.3\columnwidth}
		\centering
		\raisebox{-.25\height}{\includegraphics[width=1.2cm]{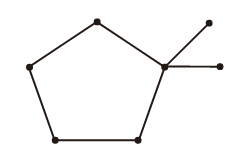}}
	\end{minipage}
    & 2.4908
    & \begin{minipage}[b]{0.3\columnwidth}
		\centering
		\raisebox{-.25\height}{\includegraphics[width=1.2cm]{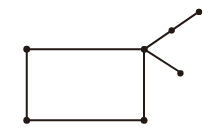}}
	\end{minipage}
    &  2.499
    \\[2.5pt]
   \hline
   \begin{minipage}[b]{0.3\columnwidth}
		\centering
		\raisebox{-.30\height}{\includegraphics[width=0.9cm]{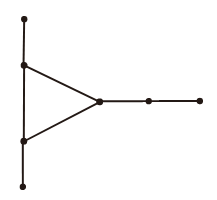}}
	\end{minipage}
    & 2.5202
    & \begin{minipage}[b]{0.3\columnwidth}
		\centering
		\raisebox{-.25\height}{\includegraphics[width=1.2cm]{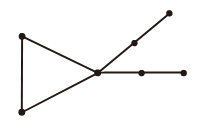}}
	\end{minipage}
    &  2.5376
    \\[2.5pt]
    \hline
     \begin{minipage}[b]{0.3\columnwidth}
		\centering
		\raisebox{-.25\height}{\includegraphics[width=1.6cm]{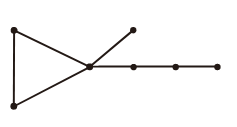}}
	\end{minipage}
    & 2.5727
    & \begin{minipage}[b]{0.3\columnwidth}
		\centering
		\raisebox{-.25\height}{\includegraphics[width=1.2cm]{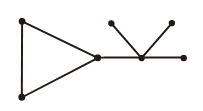}}
	\end{minipage}
    &  2.5962
    \\[2.5pt]
    \hline
    \begin{minipage}[b]{0.3\columnwidth}
		\centering
		\raisebox{-.25\height}{\includegraphics[width=1.2cm]{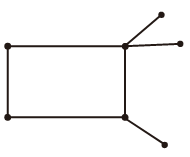}}
	\end{minipage}
    & 2.5992
    & \begin{minipage}[b]{0.3\columnwidth}
		\centering
		\raisebox{-.25\height}{\includegraphics[width=1.5cm]{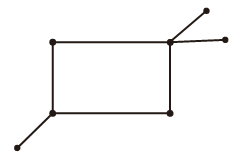}}
	\end{minipage}
    &  2.6023
    \\[2.5pt]
    \hline
     \begin{minipage}[b]{0.3\columnwidth}
		\centering
		\raisebox{-.25\height}{\includegraphics[width=1.2cm]{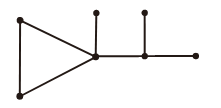}}
	\end{minipage}
    & 2.6209
    & \begin{minipage}[b]{0.3\columnwidth}
		\centering
		\raisebox{-.25\height}{\includegraphics[width=1.2cm]{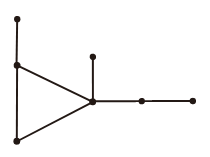}}
	\end{minipage}
    &  2.6564
    \\[2.5pt]
    \hline
     \begin{minipage}[b]{0.3\columnwidth}
		\centering
		\raisebox{-.25\height}{\includegraphics[width=1.4cm]{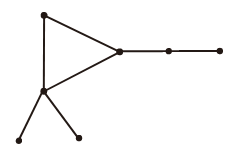}}
	\end{minipage}
    & 2.6795
    & \begin{minipage}[b]{0.3\columnwidth}
		\centering
		\raisebox{-.25\height}{\includegraphics[width=0.9cm]{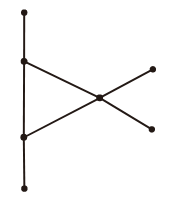}}
	\end{minipage}
    &  2.75
    \\[2.5pt]
    \hline
     \begin{minipage}[b]{0.3\columnwidth}
		\centering
		\raisebox{-.25\height}{\includegraphics[width=1.3cm]{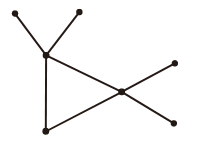}}
	\end{minipage}
    & 2.8717
    & \begin{minipage}[b]{0.3\columnwidth}
		\centering
		\raisebox{-.25\height}{\includegraphics[width=1.3cm]{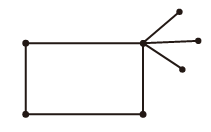}}
	\end{minipage}
    &  2.9314
    \\[2.5pt]
    \hline
     \begin{minipage}[b]{0.3\columnwidth}
		\centering
		\raisebox{-.25\height}{\includegraphics[width=1.2cm]{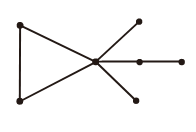}}
	\end{minipage}
    & 2.9516
    & \begin{minipage}[b]{0.3\columnwidth}
		\centering
		\raisebox{-.35\height}{\includegraphics[width=1cm]{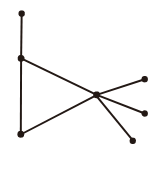}}
	\end{minipage}
    &  3.0453
    \\[2.5pt]
    \hline
    \begin{minipage}[b]{0.3\columnwidth}
		\centering
		\raisebox{-.25\height}{\includegraphics[width=1cm]{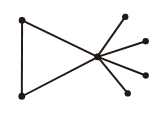}}
	\end{minipage}
    & 3.4526
    &
    &
    \\[2.5pt]
     \hline
  \end{tabular}
  \caption{The approximate value of $\rho_{ag}$ of the unicyclic graph with order $n=7$}
\end{table}

\clearpage
Suppose that $n\geq8$, the proof of $G\in\mathcal{G}_{2}$ and $G\ncong
G_{2}$, $G\in\mathcal{G}_{3}$ and $G\ncong
G_{3}$, $G\in\mathcal{G}_{4}$ and $G\ncong
G_{3}$ or $G\in\mathcal{G}_{5}$ and $G\ncong
G_{5}$ in Theorem 4.1 is given in the following.

\noindent (1) If $G\in\mathcal{G}_{2}$ and $G\ncong G_{2}$, then $d_{1}+d_{2}=n$ and $3\leq d_{1},d_{2}\leq
n-3$.

\noindent For the vertex $v_{1}$,

\begin{eqnarray*}
(A_{ag}(G)\textbf{x})_{v_{1}}&=&(d_{1}-2)\cdot\frac{1+d_{1}}{2\sqrt{d_{1}}}+\frac{d_{1}+d_{2}}{2\sqrt{d_{1}d_{2}}}\cdot\sqrt{d_{2}}+\frac{d_{1}+3}{2\sqrt{3d_{1}}}\cdot\sqrt{3}\\
&=&\frac{d_{1}^{2}+n+1}{2\sqrt{d_{1}}}\leq \frac{(n-1)d_{1}}{2\sqrt{d_{1}}}=\frac{(n-1)}{2}\cdot\sqrt{d_{1}}.
\end{eqnarray*}

\noindent Similarly, for the vertex $v_{2}$, we get
$\displaystyle (A_{ag}(G)\textbf{x})_{v_{2}}\leq \frac{(n-1)\sqrt{d_{2}}}{2}$.

\noindent For the vertex $v_{3}$,

\begin{eqnarray*}
(A_{ag}(G)\textbf{x})_{v_{3}}&=&\frac{3+d_{1}}{2\sqrt{3d_{1}}}\cdot\sqrt{d_{1}}+\frac{3+d_{2}}{2\sqrt{3d_{2}}}\cdot\sqrt{d_{2}}+\frac{1+3}{2\sqrt{3}}
=\frac{10+d_{1}+d_{2}}{2\sqrt{3}}=\frac{n+10}{2\sqrt{3}}\\
&=&\frac{n+10}{6}\cdot\sqrt{3}< \frac{n-1}{2}\cdot\sqrt{d_{3}}.
\end{eqnarray*}

Then $\displaystyle (A_{ag}(G)\textbf{x})_{i}< \frac{n-1}{2}x_{i}$ holds for all vertices in $G\in\mathcal{G}_{2}$ and $G\ncong
G_{2}$.

\noindent (2) If $G\in\mathcal{G}_{3}$ and $G\ncong G_{3}$, then $d_{1}+d_{2}=n$ and $3\leq d_{1},d_{2}\leq n-3$.

\noindent For the vertex $v_{1}$,

\begin{eqnarray*}
(A_{ag}(G)\textbf{x})_{v_{1}}&=&(d_{1}-2)\cdot\frac{1+d_{1}}{2\sqrt{d_{1}}}+\frac{d_{2}+d_{1}}{2\sqrt{d_{1}d_{2}}}\cdot\sqrt{d_{2}}+\frac{d_{1}+2}{2\sqrt{2d_{1}}}\cdot\sqrt{2}\\
&=&\frac{d_{1}^{2}+n}{2\sqrt{d_{1}}} <\frac{(n-1)d_{1}}{2\sqrt{d_{1}}}=\frac{(n-1)\sqrt{d_{1}}}{2}.
\end{eqnarray*}

\noindent For the vertex $v_{2}$,

\begin{eqnarray*}
(A_{ag}(G)\textbf{x})_{v_{2}}&=&(d_{2}-3)\cdot\frac{1+d_{2}}{2\sqrt{d_{2}}}+\frac{d_{2}+d_{1}}{2\sqrt{d_{2}d_{1}}}\cdot\sqrt{d_{1}}+2\cdot\frac{d_{2}+2}{2\sqrt{2d_{2}}}\cdot\sqrt{2}\\
&=&\frac{d_{2}^{2}+n+1}{2\sqrt{d_{2}}} \leq
\frac{(n-1)d_{2}}{2\sqrt{d_{2}}}=\frac{(n-1)\sqrt{d_{2}}}{2}.
\end{eqnarray*}

\noindent For the vertex $v_{3}$,

\begin{eqnarray*}
(A_{ag}(G)\textbf{x})_{v_{3}}&=&\frac{2+d_{1}}{2\sqrt{2d_{1}}}\cdot\sqrt{d_{1}}+\frac{2+d_{2}}{2\sqrt{2d_{2}}}\cdot\sqrt{d_{2}}\\
&=&\frac{4+d_{1}+d_{2}}{2\sqrt{2}}=\frac{4+n}{2\sqrt{2}}<
\frac{n-1}{2}\cdot\sqrt{d_{3}}.
\end{eqnarray*}

\noindent For the vertex $v_{4}$,

$\displaystyle \ \ \
(A_{ag}(G)\textbf{x})_{v_{4}}=\frac{2+d_{2}}{2\sqrt{2d_{2}}}\cdot\sqrt{d_{2}}+\frac{1+2}{2\sqrt{2}}=\frac{5+d_{2}}{2\sqrt{2}}<
\frac{n-1}{2}\cdot\sqrt{d_{4}}$.

 Then $\displaystyle (A_{ag}(G)\textbf{x})_{i}< \frac{n-1}{2}x_{i}$ holds for all vertices in $G\in\mathcal{G}_{3}$ and $G\ncong
 G_{3}$.

\noindent (3) If $G\in\mathcal{G}_{4}$ and $G\ncong G_{3}$, then $d_{1}+d_{2}=n$ and $3\leq d_{1},d_{2}\leq n-3$.

\noindent For the vertex $v_{1}$,

\begin{eqnarray*}
(A_{ag}(G)\textbf{x})_{v_{1}}&=&(d_{1}-1)\cdot\frac{1+d_{1}}{2\sqrt{d_{1}}}+\frac{d_{2}+d_{1}}{2\sqrt{d_{1}d_{2}}}\cdot\sqrt{d_{2}}\\
&=&\frac{d_{1}^{2}+n-1}{2\sqrt{d_{1}}}<
\frac{(n-1)d_{1}}{2\sqrt{d_{1}}}=\frac{(n-1)\sqrt{d_{1}}}{2}.
\end{eqnarray*}

\noindent For the vertex $v_{2}$,

\begin{eqnarray*}
(A_{ag}(G)\textbf{x})_{v_{2}}&=&(d_{2}-3)\cdot\frac{1+d_{2}}{2\sqrt{d_{2}}}+\frac{d_{2}+d_{1}}{2\sqrt{d_{1}d_{2}}}\cdot\sqrt{d_{1}}+2\cdot\frac{d_{2}+2}{2\sqrt{2d_{2}}}\cdot\sqrt{2}\\
&=&\frac{d_{2}^{2}+n+1}{2\sqrt{d_{2}}} \leq
\frac{(n-1)d_{2}}{2\sqrt{d_{2}}}=\frac{(n-1)\sqrt{d_{2}}}{2}.
\end{eqnarray*}

\noindent For the vertex $v_{3}$,

$\displaystyle
(A_{ag}(G)\textbf{x})_{v_{3}}=\frac{2+d_{2}}{2\sqrt{2d_{2}}}\cdot\sqrt{d_{2}}+\frac{2+2}{2\sqrt{4}}\cdot\sqrt{2}=\frac{d_{2}+6}{2\sqrt{2}}<
\frac{n-1}{2}\cdot\sqrt{d_{3}}$.

\noindent Similarly, for the vertex $v_{4}$, we obtain
\begin{center}
$\displaystyle (A_{ag}(G)\textbf{x})_{v_{4}}< \frac{(n-1)\sqrt{d_{4}}}{2}$.
\end{center}

Then $\displaystyle (A_{ag}(G)\textbf{x})_{i}<\frac{n-1}{2}x_{i}$ holds for all vertices in $G\in\mathcal{G}_{4}$ and $G\ncong
G_{3}$.

\noindent (4) If $G\in\mathcal{G}_{5}$ and $G\ncong G_{4}$, then $d_{1}+d_{2}=n$ and $3\leq d_{1},d_{2}\leq n-3$.

\noindent For the vertex $v_{1}$,

\begin{eqnarray*}
(A_{ag}(G)\textbf{x})_{v_{1}}&=&(d_{1}-2)\cdot\frac{1+d_{1}}{2\sqrt{d_{1}}}+\frac{d_{1}+d_{2}}{2\sqrt{d_{1}d_{2}}}\cdot\sqrt{d_{2}}+\frac{d_{1}+2}{2\sqrt{2d_{1}}}\cdot\sqrt{2}\\
&=&\frac{d_{1}^{2}+n}{2\sqrt{d_{1}}}<\frac{(n-1)d_{1}}{2\sqrt{d_{1}}}=\frac{(n-1)\sqrt{d_{1}}}{2}.
\end{eqnarray*}

\noindent Similarly, for the vertex $v_{2}$, we obtain
\begin{center}
$\displaystyle (A_{ag}(G)\textbf{x})_{v_{2}}<\frac{(n-1)\sqrt{d_{2}}}{2}$.
\end{center}

\noindent  For the vertex $v_{3}$,

\begin{center}
$\displaystyle \ \ \
(A_{ag}(G)\textbf{x})_{v_{3}}=\frac{2+d_{2}}{2\sqrt{2d_{2}}}\cdot\sqrt{d_{2}}+\frac{2+2}{2\sqrt{4}}\cdot\sqrt{2}=\frac{d_{2}+6}{2\sqrt{2}}<
\frac{n-1}{2}\cdot\sqrt{d_{3}}$.
\end{center}

\noindent Similarly, for the vertex $v_{4}$, we get
\begin{center}
$\displaystyle (A_{ag}(G)\textbf{x})_{v_{4}}<\frac{(n-1)\sqrt{d_{4}}}{2}$.
\end{center}

 Then $\displaystyle (A_{ag}(G)\textbf{x})_{i}<\frac{n-1}{2}x_{i}$ holds for all vertices in $G\in\mathcal{G}_{5}$ and $G\ncong G_{4}$.

That is $\displaystyle A_{ag}(G)\textbf{x}\leq \frac{n-1}{2}\textbf{x}$ holds for $G$.

\end{document}